\documentclass[a4paper,10pt]{article}

\usepackage{amsmath,amsthm}

\newcommand{\RR}{\mathbf{R}}
\newcommand{\norm}[1]{\|#1\|}
\newcommand{\scp}[2]{\langle #1,#2\rangle}
\newcommand{\Uad}{{U_\textup{ad}}}

\newtheorem{lemma}{Lemma}[section]
\newtheorem{proposition}[lemma]{Proposition}
\newtheorem{theorem}[lemma]{Theorem}

\theoremstyle{definition}
\newtheorem{definition}[lemma]{Definition}

\DeclareMathOperator*{\argmin}{argmin}

\author{Dirk A. Lorenz\thanks{Corresponding author, Institute for Analysis and Algebra, TU Braunschweig, D-38902 Braunschweig, Germany, \texttt{d.lorenz@tu-braunschweig.de}.} \and Arnd R\"osch\thanks{Universit\"at Duisburg-Essen, Fachbereich Mathematik, Forsthausweg 2, D-47057 Duisburg, Germany, \texttt{arnd.roesch@uni-due.de}.}}
\title{Error estimates for joint Tikhonov- and Lavrentiev-regularization of constrained control problems}

\begin{document}
\maketitle
\begin{abstract}
  We consider joint Tikhonov- and Lavrentiev-regularization of control
  problems with pointwise control- and state-constraints. We derive
  error estimates for the error which is introduced by the Tikhonov
  regularization. With the help of this results we show, that if the
  solution of the unconstrained problem has no active constraints, the
  same holds for the Tikhonov-regularized solution if the
  regularization parameter is small enough and a certain source
  condition is fulfilled.
\end{abstract}



\section{Introduction}
\label{sec:introduction}
In this paper we consider problems that can be interpreted either as
optimal control problems or identification problems (inverse
problems).  Let $D'\subset D$ be bounded domains in $\RR^N$ ($N =
2,3$), $U=L^2(D)$ and $Y$ be a Hilbert space. For a compact and linear
mapping $S:U\to Y$, an element $y_d\in Y$ and bounded measurable
functions $b:D\to \RR$ and $\psi:D'\to \RR$ we consider the
constrained minimization problem
\begin{equation}
  \label{eq:P}
  \tag{$P$}
  \min \norm{Su-y_d}^2\quad \text{ s.t. }\ 0\leq u \leq b
\text{ a.e. on }D,\quad Su\leq \psi.  \text{ a.e. on }D'
\end{equation}
The space $Y$ is called data space in inverse problems and state space for
control problems. Moreover $U$ is considered as solutions space or control space.
Our specific feature is the presence of two types of different inequality constraints.
The first one
\[
0\leq u \leq b\text{ a.e. on }D
\]
describes a set of physically meaningful solutions in inverse problems or
a set of admissible controls. The second one
\[
 Su\leq \psi.  \text{ a.e. on }D'
\]
ensures that the data (the state) is pointwise in a reasonable range.
In the inverse-problem context, the constraints on $u$ and $Su$ model
prior knowledge on the solution and the measured data, respectively
and hence, they shall lead to more reliable reconstructions.

Before we start to analyze this problem, we will give a specific example.
Let $S: u\mapsto y$ be the solution operator of the boundary value problem
\[
-\Delta y =u \text{ on } D, \quad y=0 \text{ on } \partial D
\]
with a Lipschitz domain $D$. We specify the spaces $U=Y=L^2(D)$.
The set $D'\subset D$ can be for instance an inner subdomain. However, it would also be possible
to discuss sets $D'$ containing only finitely many points. Of course our setting fits also
to more challenging problems.

Next, we will discuss the regularization of such a problem. Let us first mention that the
introduced problem~(\ref{eq:P}) can be well-posed or ill-posed. One can easily construct situations,
where the set of functions $u$ satisfying the inequality constraints is empty or consists of
exactly one point. However, such situations are not in our focus. We will assume later that  
the set of feasible functions $u\in U$ has an inner point with respect to the $L^\infty$-topology.

Since the set of feasible functions $u$ is weakly compact in $U$, the existence of at least one
function $u$ satisfying the inequality constraints ensures the existence of a solution of (\ref{eq:P}). 
Consequently, existence of solutions of (\ref{eq:P}) will not be an important issue in this paper.

Let us now think about uniqueness of solutions. Using standard arguments one can show the
uniqueness of the state (resp.~data) $y:=Su$. However, the uniqueness of $u$ is only guaranteed if $S$ is injective. 
Of course, a standard Tikhonov regularization guarantees uniqueness of solutions of the
regularized problems. Therefore, the discussion of uniqueness aspects seems to be complete.
However, there is another uniqueness aspect occurring even in the case of finite dimensional spaces.
Solution of minimization problems are analyzed and computed by means of optimality
conditions. The most convenient form of optimality conditions includes
Lagrange multipliers, i.e., dual or adjoint variables. It may happen that there exist a subset
$D''\subset D'$ where more than one inequality holds as equality. Typically this effect is connected
with nonuniqueness of the dual variables on the corresponding set $D''$. For the construction 
of specific examples we refer to \cite{GriesseMetlaRoesch09}. Let us mention that Lagrange multipliers
associated with the inequality $Su\leq \psi$ are in general only Borel measures,
see Casas \cite{casas1993stateconstraints}. This low regularity of the dual variables was the motivation in
\cite{meyer2005stateconstraints} to introduce a second Lavrentiev type regularization.

Motivated by this argumentation we will study in this paper a family of
Tikhonov and Lavrentiev regularized problems:
\begin{equation}
  \label{eq:P_alpha_lambda}
  \tag{$P_\alpha^\lambda$}
  \min \norm{Su-y_d}^2 + \alpha\norm{u}^2\ 
\text{ s.t. } 0\leq u\leq b, \text{ a.e. on }D,
\lambda u + Su\leq \psi \text{ a.e. on }D'
\end{equation}
with $\alpha,\lambda\ge 0$.  
We denote
with $\bar u$ a solution of~(\ref{eq:P}) and with $u_\alpha^\lambda$ a
solution of~(\ref{eq:P_alpha_lambda}). Note that the solution
of~(\ref{eq:P_alpha_lambda}) is always unique if it exists.

Next, one can study a lot of different regularization errors. People coming from inverse
problems usually study the solution behavior for $\alpha\downarrow 0$. Since the 
inverse theory heavily uses representation formulas, pointwise inequalities are undesired. 

People coming from optimal control state that the Tikhonov term
represents the costs for the control and assume that the parameter
$\alpha$ is a given, fixed model parameter. Hence, it would be enough
to study the behavior for $\lambda\downarrow 0$ for a fixed Tikhonov
parameter $\alpha>0$. However, even in the optimal-control context one
may take the position that the problem~(\ref{eq:P}) is the one which
shall be solved and that the Tikhonov regularization is employed only
to stabilize the problem.

Of course, error estimates for the whole regularization process  are of high interest.
However, this problem seems to be challenging. In this paper, 
we will only contribute a little bit in answering this question.

Let us summarize our plans.
The  problem~(\ref{eq:P}) is the one which we want to
solve.  Therefore, we are interested in estimates for the error $\bar u
- u_\alpha^\lambda$.  We split this error into the \emph{Tikhonov
  error} and the \emph{Lavrentiev error}:
\[
\norm{\bar u - u_\alpha^\lambda} \leq \underbrace{\norm{\bar u -
    u_\alpha^0}}_{\text{Tikhonov error}} + \underbrace{\norm{u_\alpha^0 -
    u_\alpha^\lambda}}_\text{Lavrentiev error}.
\]
Estimates for the Lavrentiev error can be found in
\cite{meyer2005stateconstraints,cherednichenko2008lavrentieverror,cherednichenko2008lavrentiev,hinze2008lavrentiev}.
Here, we focus on estimates of the Tikhonov error.  The analysis of
this error is well developed in the framework on inverse problems for problems
without inequality constraints.
However, for the constrained case there are only few
results~\cite[Section 5.4]{engl1996inverseproblems}. 


The paper is structured as follows: In Section 2 we will state some preliminary results.
Error estimates for the constrained Tikhonov regularization are located in Section 3.
The activity of inequality constraints is analyzed in Section 4. Error estimates
for the Lavrentiev regularization can be found in Section 5. The verification of the
assumptions for a distributed elliptic control problem and a Fredholm integral operator is presented in Section 6.

\section{Preliminary results}
\label{sec:preliminary-results}
First, we introduce  notations for the admissible sets for~(\ref{eq:P})
and~(\ref{eq:P_alpha_lambda}) 
\begin{align*}
  \Uad & = \{u\in U\ |\ 0\leq u\leq b,\ Su\leq \psi\},\\
  \Uad^\lambda & = \{u\in U\ |\ 0\leq u\leq b,\ \lambda u+Su\leq \psi\}.
\end{align*}
Note that, due to $0\leq u$ we have
$\Uad^\lambda\subset \Uad^\mu$ for $\lambda\geq \mu$.  With this
notation we reformulate
\begin{align}
  \label{eq:Pvar} 
  \tag{$P$} 
  \min \norm {Su-y_d}^2\ \text{ s.t. }\ u\in \Uad\\
  \label{eq:P_alpha_lambda_var}
   \tag{$P_\alpha^\lambda$}
   \min \norm{S u-y_d}^2 + \alpha\norm{u}^2\ \text{ s.t. }\ u\in \Uad^\lambda.  
\end{align} 

Now we state some preliminary results:
\begin{lemma}
  Let the operator $S$ be linear and continuous.
  Then the sets $\Uad$ and $\Uad^\lambda$ are closed, convex and bounded.
\end{lemma}

\begin{definition}
  The \emph{constrained pseudo inverse} of an operator $S$ with
  respect to a convex and closed set $C$ is defined via
  \[\norm{S_C^+(y_d)} = \min\{\norm{u}\ |\ u\in\argmin_C \norm{Su-y_d}^2\}.\]
  In other words, $S_C^+(y_d)$ is the minimizing element of the
  residuum which has minimal norm.
\end{definition}
The operator $S_C^+$ is a nonlinear operator with the following
properties.
\begin{lemma}
  If $C$ is non-empty, convex, closed and bounded, then
  \[D(S_C^+) = Y.\]
\end{lemma}
\begin{proof}
  Note that $C$ is weakly sequentially closed and hence a minimizing
  sequence of $\norm{Su-y_d}^2$ in $U$ has a weak accumulation point.
  By lower semicontinuity we see that this accumulation point is
  indeed a minimizer, and hence, the constrained pseudo inverse
  exists.
\end{proof}
The following well known proposition shows continuity in the Tikhonov
parameter $\alpha$:
\begin{proposition}\label{prop:cont_dep_alpha}
$u_\alpha^\lambda$ depends continuously on $\alpha$ for $\alpha>0$.  
\end{proposition}
\begin{proof}
  We drop the superscript $\lambda$ because it is fixed here.
  Consider $u_\alpha$ and $u_\beta$ for $\alpha,\beta>0$.  Since
  $\lambda$ is fixed, $u_\beta$ is admissible for
  (\ref{eq:P_alpha_lambda_var}) and vice versa.  Hence, we can insert
  these elements in the corresponding variational inequalities and
  obtain
  \begin{gather*}
    \scp{S^*Su_\alpha + \alpha u_\alpha -
      S^*y_d}{u_\beta-u_\alpha}\geq 0\\
    \scp{S^*Su_\beta + \beta u_\beta -
      S^*y_d}{u_\alpha-u_\beta}\geq 0.
  \end{gather*}
  We get
  \[
  \scp{S^*Su_\alpha + \alpha u_\alpha - S^*Su_\beta - \beta
    u_\beta}{u_\beta-u_\alpha}\geq 0
  \]
  and obtain
  \[
  \norm{Su_\alpha - Su_\beta}^2\leq \scp{\alpha u_\alpha - \beta
    u_\beta}{u_\beta-u_\alpha}
  \leq |\alpha-\beta| \norm{u_\alpha}\norm{u_ \beta-u_\alpha} - \beta\norm{u_\beta-u_\alpha}^2.
  \]
  Since the left hand side is positive we get
  \[
  \norm{u_\beta-u_\alpha} \leq \frac{|\alpha-\beta|}{\beta}\norm{u_\alpha}
  \]
  and the right hand side converges to 0 for $\beta\to\alpha$.
\end{proof}
Finally, we state a characterization of the solution
of~(\ref{eq:P_alpha_lambda_var}) by means of a projection formula:
\begin{lemma}[Projection formula]
  The solution $u_\alpha^\lambda$ of~(\ref{eq:P_alpha_lambda_var}) is
  characterized by
  \begin{equation}
    \label{eq:projection_formula}
    u_\alpha^\lambda = P_{\Uad^\lambda}\Bigl(-\frac{S^*(Su_{\alpha}^\lambda - y_d)}{\alpha}\Bigr).
  \end{equation}
\end{lemma}
\begin{proof}
  We rewrite~(\ref{eq:P_alpha_lambda_var}) with the help of the
  indicator function $I_{\Uad^\lambda}$ as
  \[
  \min_u \norm{Su-y_d}^2 + \alpha\norm{u}^2 + I_{\Uad^\lambda}(u).
  \]
  With the help of subgradient calculus we get from optimality of
  $u_\alpha^\lambda$
  \[
  -\frac{S^*(Su_\alpha^\lambda - y_d)}{\alpha} \in
  \partial(\norm{\cdot}^2+I_{\Uad^\lambda})(u_\alpha^\lambda).
  \]
  Since $\partial(\norm{\cdot}^2+I_{\Uad^\lambda})^{-1} =
  P_{\Uad^\lambda}$ (see, e.g.~\cite{rockafellar1998variation}) this
  gives
  \[
  P_{\Uad^\lambda}\Bigl(-\frac{S^*(Su_\alpha^\lambda -
    y_d)}{\alpha}\Bigr) = u_\alpha^\lambda.
  \]
\end{proof}

\section{Error estimates for constrained Tikhonov regularization}
\label{sec:error-estim-constr}
Now we establish error estimates for the term $\norm{\bar u -
  u_\alpha^0}$.  Resembling results can be found
in~\cite{engl1996inverseproblems}. Due to the structural differences
between inverse problems and optimal control we state the result in
our terminology with an explicit constant and present a full proof.
\begin{theorem}
  \label{thm:conv_rate_alpha}
  Let $\bar u$ be a solution of (\ref{eq:P}) and denote with $P_{\Uad}$ the projection onto $\Uad$.
  Moreover, let the following source
  condition be fulfilled:
  \[
  \exists w\in Y:\ \bar u = P_{\Uad}(S^* w).
  \]
  Then it holds:
  \begin{eqnarray}
    \label{eq:err_estimate_u_alpha_0}
    \norm{u_\alpha^0 - \bar u} & \leq &\sqrt{\alpha}\norm{w} + \frac{1}{\sqrt{\alpha}}\norm{S\bar u - y_d}\\
    \label{eq:err_estimate_Su_alpha_0}
    \norm{Su_\alpha^0 - y_d}  &\leq &2\alpha\norm{w} + \norm{S \bar u - y_d}    
  \end{eqnarray}
\end{theorem}
\begin{proof}
  By definition we have $u_\alpha^0\in \Uad$ and $\bar u\in \Uad$.
  Hence, by optimality of $u_\alpha^0$ we have
  \[
  \norm{Su_\alpha^0 - y_d}^2 + \alpha\norm{u_\alpha^0}^2
  \leq \norm{S\bar u - y_d}^2 + \alpha\norm{\bar u}^2.
  \]
  Rearranging yields
  \begin{equation*}
    \norm{Su_\alpha^0 - y_d}^2 + \alpha(\norm{u_\alpha^0}^2 - \norm{\bar u}^2)\leq \norm{S\bar u - y_d}^2
  \end{equation*}
  which we extend to
  \begin{equation}
    \label{eq:aux_est1}
    \norm{Su_\alpha^0 - y_d}^2 + \alpha(\norm{u_\alpha^0}^2 - \norm{\bar u}^2 - 2\scp{S^*w}{u_\alpha^0 - \bar u} + 2\scp{S^*w}{u_\alpha^0 - \bar u})
    \leq \norm{S\bar u - y_d}^2.
  \end{equation}
  Since $\bar u = P_\Uad S^*w$ we have for all $u\in\Uad$
  \[
  \scp{S^*w - \bar u}{u-\bar u}\leq 0.
  \]
  Using this with $u=u_\alpha^0$ we get
  \[
  \norm{u_\alpha^0}^2 - \norm{\bar u}^2 - 2\scp{S^*w}{u_\alpha^0 -
    \bar u} \geq \norm{u_\alpha^0}^2 - \norm{\bar u}^2 - 2\scp{\bar
    u}{u_\alpha^0 - \bar u} = \norm{u_\alpha^0-\bar u}^2.
  \]
  We further estimate from~(\ref{eq:aux_est1})
  \begin{align*}
    & & \norm{Su_\alpha^0 - y_d}^2 + \alpha(\norm{u_\alpha^0 - \bar u}^2  + 2\scp{S^*w}{u_\alpha^0 - \bar u})
    & \leq \norm{S\bar u - y_d}^2\\
    & \Leftrightarrow & \norm{Su_\alpha^0 - y_d}^2 + 2\scp{\alpha w}{Su_\alpha^0 - y_d} + \alpha\norm{u_\alpha^0 - \bar u}^2
    & \leq \norm{S\bar u - y_d}^2 + 2\scp{\alpha w}{S\bar u -y_d}.
  \end{align*}
  Completing the squares by adding $\norm{\alpha w}^2$ gives
  \[
  \norm{Su_\alpha^0 - y_d + \alpha w}^2 + \alpha\norm{u_\alpha^0 - \bar u}^2
  \leq \norm{S \bar u - y_d + \alpha w}^2.
  \]
  On the one hand this leads to
  \[
  \norm{u_\alpha^0 - \bar u}^2
  \leq  \tfrac{1}{\alpha}\norm{S \bar u - y_d +\alpha w}^2
  \]
  which gives by taking square roots and using the triangle inequality
  \[
  \norm{u_\alpha^0 - \bar u} \leq \tfrac{1}{\sqrt{\alpha}}\norm{S \bar
    u - y_d} + \sqrt{\alpha}\norm{w}.
  \]
  On the other hand we conclude
  \[
  \norm{Su_\alpha^0 - y_d + \alpha w}^2 
  \leq \norm{S \bar u - y_d + \alpha w}^2
  \]
  which implies
  \[
  \norm{Su_\alpha^0 - y_d}
  \leq \norm{S \bar u - y_d}  + 2\alpha\norm{w}.  
  \]
\end{proof}
The estimate~(\ref{eq:err_estimate_u_alpha_0}) motivates the following
parameter choice: If $y_d$ is not in the range, we see that the right
hand side in~(\ref{eq:err_estimate_u_alpha_0}) is smallest for
\[
\alpha^* = \frac{\norm{S \bar u - y_d}}{\norm{w}}
\]
and hence, is a reasonable choice for the choice of the Tikhonov
parameter $\alpha$ if the quantities were known.  Nonetheless, the
error estimate~(\ref{eq:err_estimate_u_alpha_0}) is useful for the
determination of the total error $\norm{u_\alpha^\lambda - \bar u}$.

\section{Activity of the constraints}
\label{sec:activity-constraints}
In this section we investigate the following situation: Assume that
the optimal solution $\bar u$ of~(\ref{eq:Pvar}) has no active inequality
constraints. That means that in fact we would get the same solution
without imposing any inequality constraints. However, the formulation with
additional inequality constraints is reasonable since a solution of an unconstrained
inverse problem can violate these constraints for noisy data.

Now the question arises: Can we expect a solution
without active constraints also for the purely Tikhonov-regularized
problem for small regularization parameters? In fact this can be shown
with the help of the estimates of Theorem~\ref{thm:conv_rate_alpha}:
\begin{theorem}
  \label{thm:u_alpha_in_intUad}
  Let $S:U=L^2(D)\to Y$ with dense range and assume that
  \begin{align}
    \label{eq:ass_S}
    u_n\to u \text{ in } U & \implies Su_n\to Su \text{ in } L^\infty(D')\\
    \label{eq:ass_S*}
    y_n \rightharpoonup y \text{ in } Y & \implies S^*y_n \to Sy \text{ in } L^\infty(D)
  \end{align}
  Let $\bar u$ be a
  solution of~(\ref{eq:Pvar}) such that for some $\tau>0$ it holds
  that $\tau < \bar u <b-\tau$ and $S\bar u<\psi-\tau$ hold.
  Moreover, let there be $w\in Y$ such that $\bar u = S^* w$.
  
  Then there exists $\alpha_0>0$ such that for every $\alpha<\alpha_0$
  the solution $u_\alpha^0$ of~($P_\alpha^0$) also fulfills $\tau <
  u_\alpha^0 < b-\tau$ and $Su_\alpha^0<\psi-\tau$.
\end{theorem}
\begin{proof}
  Since $\bar u$ does not have active constraints and $S$ has dense
  range, it holds $S\bar u = y_d$.  To see this, note that $\bar u$
  fulfills
  \[
  \scp{S^*(S\bar u-y_d)}{u-\bar u}\geq 0\quad \text{for all}\quad u\in\Uad
  \]
  and let us assume that there exists $v\in U$ but $v\notin \Uad$, such that
  \begin{equation}
    \label{eq:var_ineq_v}
    \scp{S^*(S\bar u-y_d)}{v-\bar u}<0.
  \end{equation}
  We may approximate $v$ by bounded functions and hence, there exists
  a bounded $\bar v$ such that
  \[
  \scp{S^*(S\bar u-y_d)}{\bar v-\bar u}<0.
  \]
  Now observe that for $\theta > 0$ small enough, the function
  $v^\theta = \theta \bar v +(1-\theta)\bar u$ is in $\Uad$ since
  $\bar u$ does not have active constraints. We see
  \[
  \scp{S^*(S\bar u-y_d)}{v^\theta-\bar u} = \theta\scp{S^*(S\bar
    u-y_d)}{\bar v-\bar u}<0
  \]
  which contradicts the optimality of $\bar u$
  for~\eqref{eq:Pvar}. Hence, \eqref{eq:var_ineq_v} has to be
  fulfilled for all $v\in U$ and this shows that $\bar u$ is also a
  solution of the unconstrained problem. Since $S$ has dense range,
  this optimal value for this problem is 0 which shows $S\bar u =
  y_d$.

  Moreover, we have by assumption $\bar u = S^* w = P_\Uad
  S^*w$.  We conclude from Theorem~\ref{thm:conv_rate_alpha}
  that
  \begin{align*}
    \norm{u_\alpha^0 - \bar u} & \leq \sqrt{\alpha}\norm{w}.
  \end{align*}
  This implies $u_\alpha^0\to \bar u$ in $L^2$ and by assumption we know
  that $Su_\alpha^0 \to S\bar u$ in $L^\infty(D')$.  Because of
$S\bar u<\psi-\tau$, we find $Su_\alpha^0<\psi-\tau/2$ for sufficiently small $\alpha$. 
  Hence, we can use the
  projection formula~(\ref{eq:projection_formula}) with $P_{[0,b]}$
  instead of $P_{\Uad}$ for small $\alpha$ and get
  \[
  u_\alpha^0 = P_{[0,b]}(-S^*\tfrac{1}{\alpha}(Su_\alpha^0 - y_d)).
  \]
  Again from Theorem~\ref{thm:conv_rate_alpha} we find
  \[
  \norm{Su_\alpha^0 - y_d}  \leq 2\alpha\norm{w}.
  \]
  Consequently,
  $\norm{(Su_\alpha^0 - y_d)/\alpha}$ is uniformly bounded in $Y$.
  Let us now take an arbitrary sequence $\{\alpha_n\}$ with $\alpha_n\to 0$ for
  $n$ to $\infty$. Since $(Su_{\alpha_n}^0 - y_d)/\alpha_n$ is uniformly
  bounded in $Y$,
  we can find a weakly convergent subsequence in $Y$ denoted by the index $n'$.
  The assumption on $S^*$ yields
  strong convergence of  $u_{\alpha_{n'}}^0$ in $L^\infty(D)$.
  We already know that $u_{\alpha_{n'}}^0$ converges in $U=L^2(D)$
  with limit $\bar u$. By a standard argumentation we obtain that
  $u_\alpha^0\to \bar u$ in $L^\infty(D)$ for $\alpha\to 0$.
  Consequently, the control constraints are also inactive for sufficiently small
  $\alpha$.
\end{proof}

Let us note that the first argumentation in the proof works also for the
weaker assumption that $u_n\to u$ in $L^p(D)$ implies $Su_n\to Su$ in $L^\infty(D')$
for a sufficiently large $p$. However, the assumptions for the adjoint operator
cannot be weakened. Thus, the practical benefit of that generalization is only small.

The above result is also of interest in the case of ill-posed
problems.  Here we may ask the question if additional constraints on
the solution in the minimization of the Tikhonov functional may
destroy the optimal convergence rate. As we will see below, this is
not the case in our setting.
Similar to Theorem~\ref{thm:conv_rate_alpha} we can state the following result:
\begin{theorem}
  Let $\bar u$ be a solution of (\ref{eq:P}) with $S\bar u = y_d$ and
  let $y^\delta$ be such that $\norm{y^\delta-y_d}\leq \delta$.
  Moreover, let the following source condition be fulfilled:
  \[
  \exists w\in Y:\ \bar u = P_{\Uad}(S^* w).
  \]
  Then it holds for
  \begin{equation}
    \label{eq:def_u_alpha_delta}
    u_\alpha^\delta = \argmin_{u\in\Uad}\norm{Su-y^\delta}^2 + \alpha\norm{u}^2
  \end{equation}
  that
  \begin{eqnarray*}
    \label{eq:err_estimate_u_alpha_delta}
    \norm{u_\alpha^\delta - \bar u} & \leq &\sqrt{\alpha}\norm{w} + \frac{\delta}{\sqrt{\alpha}}\\
    \label{eq:err_estimate_Su_alpha_delta}
    \norm{Su_\alpha^\delta - y_d}  &\leq &2\alpha\norm{w} + \delta
  \end{eqnarray*}
\end{theorem}
It also holds that all constraints get inactive for small $\alpha$ if
the true solution of noisy data does not have active constraints:
\begin{theorem}
  \label{thm:u_alpha_delta_in_intUad}
  Let $S:U\to Y$ fulfill the conditions (\ref{eq:ass_S})
  and~(\ref{eq:ass_S*}).  Let $\bar u$ be a solution
  of~(\ref{eq:Pvar}) such that for some $\tau>0$ it holds that $\tau <
  \bar u <b-\tau$ and $S\bar u<\psi-\tau$ holds. Moreover, let there
  be $w\in Y$ such that $\bar u = S^* w$ and let $u_\alpha^\delta$ be
  defined by~(\ref{eq:def_u_alpha_delta}).  Finally let
  $\alpha(\delta)$ be a parameter choice rule such that
  \[
  \alpha(\delta) \to 0,\quad \frac{\delta}{\alpha(\delta)}\to 0\ 
  \text{ for }\ \delta\to 0.
  \]
  Then there exists $\delta_0>0$ such that for $\delta<\delta_0$ it
  holds that $\tau < u_{\alpha(\delta)}^\delta <b-\tau$ and
  $Su_{\alpha(\delta)}^\delta<\psi-\tau$ holds.
\end{theorem}
\begin{proof}
  Due to the parameter choice we get $u_{\alpha(\delta)}^\delta\to
  \bar u$ strongly in $L^2$ for $\delta\to 0$. Similar to the proof of
  Theorem~\ref{thm:u_alpha_in_intUad} we conclude that the state
  constraints are not active for sufficiently small $\delta$ and hence,
  the projection formula
  \[
  u_{\alpha(\delta)}^\delta =
  P_{[0,b]}(-S^*(\tfrac{1}{\alpha(\delta)}(Su_{\alpha(\delta)}^\delta
  - y^\delta))
  \]
  holds.   Now the claim follows similarly to
  Theorem~\ref{thm:u_alpha_in_intUad}.
\end{proof}

Finally we state the following converse result which shows that the
source condition $\bar u = P_\Uad S^* w$ follows from weak convergence
of the regularized solutions together with a mild decay of the
discrepancy.
\begin{theorem}
  Let $u_\alpha^0 \rightharpoonup \bar u$ and $\norm{Su_\alpha^0 -y_d}
  = O(\alpha)$ for $\alpha\to 0$. Moreover, we assume that $S$
  fulfills~(\ref{eq:ass_S}) and (\ref{eq:ass_S*}). Then there is a
  function $w$ such that $P_\Uad S^*w = \bar u$.
\end{theorem}
\begin{proof}
  Since $\norm{Su_\alpha^0 -y_d} = O(\alpha)$ we have that
  \[
  \limsup_{\alpha\to 0} \frac{\norm{Su_\alpha^0 -y_d}}{\alpha} < \infty.
  \]
  Hence there is a sequence $\alpha_n\to 0$ and an element $w$ such that
  \[
  \frac{Su_{\alpha_n} -y_d}{\alpha_n} \rightharpoonup -w.
  \]
  By the projection formula~(\ref{eq:projection_formula}) we have
  \[
  P_\Uad\Bigl(-\frac{S^*(Su_{\alpha_n}^0 - y_d)}{\alpha_n}\Bigr) = u_{\alpha_n}^0.
  \]
For the right hand side converges weakly by assumption. With a discussion similar to that one
in Theorem 4.1 we obtain the strong convergence of the left hand side of the last equation. 
Since the both sides converge weakly  we have by uniqueness of the weak limit
  \[
  P_\Uad S^* w = \bar u.
  \]
\end{proof}

\section{Lavrentiev regularization}
We discuss now the additional Lavrentiev regularization. The motivation of this
second regularization is to improve the properties of the adjoint problem.
In this section we will investigate two different situations. In the first part we will sketch the
estimation of the Lavrentiev error in the general case. In the second one we will
again discuss the situation of Theorem~\ref{thm:u_alpha_in_intUad}. Then we will
be able to show better convergence results.

Let us start with the general case. The discussion for fixed $\alpha>0$ and a specific problem
can be found in \cite{cherednichenko2008lavrentiev}. Next, we reformulate the assumptions of
that paper for our more general setting.
\begin{gather}
  \label{eq:ass_S_var}
  S:U=L^2(D)\to L^\infty(D') \text{ continuously}.\\
  \label{eq:ass_Slater_point}
  \text{There exists } \hat u\in U, \tau>0 \text{ such that } 0\leq \hat u \leq b,\ S\hat u \leq \psi-\tau.
\end{gather}
The assumption~(\ref{eq:ass_Slater_point}) means that there exists a
Slater point with respect to the state constraints.  In the second
part of this section we will benefit from the stronger assumption that
$\bar u$ itself has this Slater property. In Section 2 we already
mentioned that the set $\Uad^\lambda$ of admissible $u$ becomes
smaller when $\lambda$ becomes larger. However, one can show the
existence of at least one admissible control for $\lambda\le
\frac{\tau}{\|\hat u\|_{L^\infty(D')}}$.

The error estimates are obtained by the following technique:
\begin{enumerate}
 \item Write down the necessary optimality conditions for $u_\alpha^0$ and $u_\alpha^\lambda$
as variational inequalities.
\item Take $u_\alpha^\lambda$ as test function in the optimality condition of $u_\alpha^0$.
\item Construct a convex linear combination $u_\sigma=\sigma \hat u+ (1-\sigma)u_\alpha^0$
which belongs to $\Uad^\lambda$. Take this function as test function in the optimality condition of $u_\alpha^\lambda$.
\item Add both inequalities and estimate all terms.
\end{enumerate}
The resulting error estimate can be found in \cite[Theorem
5.4]{cherednichenko2008lavrentiev}:
\begin{lemma} \label{lemma:lavr}
Let the assumptions~(\ref{eq:ass_S_var}) and (\ref{eq:ass_Slater_point}) be fulfilled.
Then there exists a constant $c>0$ such that for $\lambda\le \frac{\tau}{\|\hat u\|_{L^\infty(D')}}$ the error of the Lavrentiev regularization can be estimated by
\[
\norm{u_\alpha^0 -
    u_\alpha^\lambda}\le c\frac{\lambda}{\alpha}.
\]
\end{lemma}
This general result has an essential drawback: If $\alpha$ becomes small, then $\lambda$ has
to be very small to ensure a certain accuracy. 

Let us now assume that $\bar u$ itself has the Slater property
\[
S\bar u<\psi-\tau.
\]
In contrast to Section 4 we do not require an inner point property with respect to the control
constraints.

\begin{theorem} \label{thm:u_alpha_lambda}
Let $S$ fulfill assumption~(\ref{eq:ass_S_var}) and
let $\bar u$ be a solution of~(\ref{eq:P}) such that the Slater condition $ S\bar u<\psi-\tau$ and a source condition
$\bar u = P_{\Uad}(S^* w)$ are fulfilled. Then it holds for sufficiently small $\lambda$ and $\alpha$ that
\[
u_\alpha^0 = u_\alpha^\lambda.
\]
\end{theorem}
\begin{proof} From Theorem~\ref{thm:conv_rate_alpha} we know the convergence of
$\bar u_\alpha^0$ to $\bar u$ in $U$. Using the mapping property $S:\,U=L^2(D)\to
L^\infty(D')$, we can show as in the proof of
of Theorem~\ref{thm:u_alpha_in_intUad}
that
\[
S u_\alpha^0<\psi-\tau/2.
\]
holds for sufficiently small $\alpha$. From Lemma~\ref{lemma:lavr} we know that
$u_\alpha^\lambda$ tends to $u_\alpha^0$ for $\lambda\to 0$. Due to the properties of
$S$ we find
\[
S u_\alpha^\lambda<\psi
\]
for sufficiently small $\lambda$. 
For $\lambda\le \frac{\tau}{2b}$ we get
\[
S u_\alpha^0+\lambda u_\alpha^0 <\psi.
\]
Consequently, $u_\alpha^0$ is feasible for the problem ($P_\alpha^\lambda$) and 
$u_\alpha^\lambda$ is feasible for the problem ($P_\alpha^0$) for $\alpha$ and $\lambda$ small enough. Testing the optimality condition for ($P_\alpha^\lambda$) and ($P_\alpha^0$) with $u_\alpha^0$ and $u_\alpha^\lambda$, respectively, yields, similar to the proof of Proposition~\ref{prop:cont_dep_alpha},
that $u_\alpha^\lambda=u_\alpha^0$.
\end{proof}

Let us remark that the Lavrentiev regularization is also used with different sign, i.e.,
\[
Su-\lambda u\le \psi.
\]
Then, the set of admissible controls becomes larger for larger $\lambda$. Hence, the
Slater condition ensures the existence of feasible controls for arbitrary positive $\lambda$.
In this case we have no smallness condition for $\lambda$. The general result for the
Lavrentiev regularization can be found in \cite[Theorem 3.3]{cherednichenko2008lavrentiev}.
The discussion for the specific case of Theorem~\ref{thm:u_alpha_lambda} can be done
completely analogue.

Both approaches have their specific advantages. As we have already seen, the first approach
(plus sign)
yields solutions that are feasible for the unregularized minimization problem. However, we have
to deal with a smallness condition for the Lavrentiev parameter $\lambda$. The second approach
(minus sign)
does not need an additional smallness condition, but the solutions are in general not
feasible. 

Finally we combine the results of Theorem~\ref{thm:conv_rate_alpha},
Theorem~\ref{thm:u_alpha_in_intUad} and
Theorem~\ref{thm:u_alpha_lambda} to obtain an estimate for the total
error of joint Tikhonov-Lavrentiev regularization under strong
assumptions:
\begin{theorem}
  Let $S$ fulfills the assumptions~(\ref{eq:ass_S}) and (\ref{eq:ass_S*}) and
  let $\bar u$ be a solution of (\ref{eq:P}) such that there exists
  $\tau>0$ such that $\tau<\bar u<b-\tau$, $S\bar u<\psi-\tau$ and a
  $w\in Y$ such that $\bar u = S^*w$. Then it holds for $\lambda$
  small enough that the solution $u_\alpha^\lambda$
  of~(\ref{eq:P_alpha_lambda}) fulfills
  \[
  \norm{\bar u - u_\alpha^\lambda} =
  \mathcal{O}(\sqrt{\alpha})\quad\text{for}\quad\alpha\to 0.
  \]
\end{theorem}
\begin{proof}
  We split the total error as
  \[
  \norm{\bar u - u_\alpha^\lambda}\leq \norm{\bar u-u_\alpha^0}+\norm{u_\alpha^0-u_\alpha^\lambda}
  \]
  and observe that due to Theorem~\ref{thm:u_alpha_lambda} the second
  term vanishes if $\lambda$ is small enough and that the first term
  behaves like $\mathcal{O}(\sqrt{\alpha})$ due to
  Theorem~\ref{thm:conv_rate_alpha} and
  Theorem~\ref{thm:u_alpha_in_intUad}.
\end{proof}

\section{Verification of the assumption in special cases}
\subsection{An elliptic control problem}
In this section we will discuss the example from the introduction:
Let $S: u\mapsto y$ be the solution operator of the boundary value problem
\[
-\Delta y =u \text{ on } D, \quad y=0 \text{ on } \partial D
\]
with a Lipschitz domain $D\subset \RR^N$, $N\in\{2,3\}$. We specify the spaces $U=Y=L^2(D)$.
The set $D'\subset D$ is assumed to be an inner subdomain. Let us now check
the assumptions:
\begin{enumerate}
\item By the Lax-Milgram Lemma we obtain easily the existence of a
  unique solution $y\in H^1_0(D)$ for every $u\in U$. Due to the
  embedding $H^1_0(D) \hookrightarrow Y$, the operator $S$ is well
  defined

\item Let us first mention that the operator $S$ is selfadjoint for
  our specific choice of spaces.  For the uniform boundedness of
  solutions of the elliptic problem we refer to Stampacchia
  \cite{sta65}.

\item The mapping property in Theorem~\ref{thm:u_alpha_in_intUad} that
  $y_n\rightharpoonup y$ in $Y$ implies $S^*y_n \to S^*y$ in
  $L^\infty(D)$ can be obtained by the following argumentation.  Weak
  convergence $y_n\rightharpoonup y$ in $Y=L^2(D)$ implies strong
  convergence in $y_n\to y$ in $W^{-1,p}(D)$ for $2\le p< \infty$ for
  $N=2$ and $2\le p\le 6$ for $N=3$. Now the desired result follows
  again from \cite[Theorem~4.2]{sta65}.

\item The operator $S$ has dense range in $Y$. That can be verified by
  the following argumentation.  The space $C^\infty_0(D)$ is dense in
  $Y$. Moreover, all these functions belong to the image of $S$.

\item The Slater condition can only be checked (analytically or
  numerically) if one specifies the data. The Slater property for
  $\bar u$ is like the source condition an a priori assumption for the
  solution.
\end{enumerate}

Consequently, we can apply all results of our paper to that example. Only the
Slater property and the source condition are a priori assumptions. All other
assumptions of that paper are satisfied for our example.

\subsection{A Fredholm integral operator}
Another class of examples in which the properties 1.--4. are fulfilled
and which models several inverse problems is given by Fredholm
integral operators~\cite{groetsch84} or \cite[Chapter
VI.]{dunford58}. Let $D$ be a bounded domain in $\RR^N$ and let
$U=Y=L^2(D)$. For a Lipschitz continuous function $k:D\times D\to \RR$
we consider $Su(x) = \int_D k(x,x')u(x') dx'$. Let us check the
assumptions:
\begin{enumerate}
\item Since $D$ is bounded, continuity of $k$ implies that
  $S:L^2(D)\to L^2(D)$ compactly.

\item The uniform convergence of $Su_n$ to $Su$ follows from $u_n\to
  u$ by
  \begin{align*}
    \norm{Su_n-Su}_{\infty} & \leq \sup_{x\in D}\int_D |k(x,x')||u_n(x')-u(x')|dx'\\
    & \leq \sup_{x,y\in D}|k(x,y)| \sqrt{|D|}\norm{u_n-u}_2.
  \end{align*}

\item Uniform convergence of $S^*y_n$ to $S^*y$ follows from
  $y_n\rightharpoonup y$ by observing that $\norm{y_n-y}_2$ is bounded
  and that the mapping $x\mapsto \int_D k(x',x)y(x')dx'$ is Lipschitz
  continuous for every $y \in L^2(D)$.

\item Since $S$ is compact is has a representation via its singular
  value decomposition $Su = \sum_n \sigma_n\scp{u}{\psi_n}\phi_n$ with
  non-negative singular values which decay to zero. We see that the
  range of $S$ is dense if $S$ does not have zero as a singular value
  (i.e.~it is injective).

\end{enumerate}

Again, the Slater condition and the source condition are a priori
assumption for the solution.

\section*{Conclusion}
In this paper we studied the simultaneous Tikhonov and
Lavrentiev regularization of an optimal control problem with control
and state constraints. We derived error estimates in the general case
and the main tool was a source condition which resembles the classical
one used in the inverse-problem context.  With the help of this error
estimate we could prove, under additional assumptions, that for
sufficiently small Tikhonov parameter the Tikhonov-regularized
solution does not have active inequality constraint if the original
solution has the same property. Moreover, it was shown that for the
Tikhonov- and Lavrentiev-regularized problem the solution coincides
with the Tikhonov-regularized solution for sufficiently small
regularization parameters if the unregularized solution has a Slater
property. One may conclude that using additional physically motivated
inequality constraints in the context of the regularization of inverse
problems is a good idea for two reasons: For larger regularization
parameters they may yield reconstructions that are more meaningful and
asymptotically they do not destroy optimal convergence rates (and even
become inactive if they will be in the limit). On the other hand,
additional inequality constraints lead to minimization problem which
may be harder to solve. However, there are powerful algorithms
available which allow a numerical treatment of such problems
(cf.~\cite{hintermueller2003primaldualssn,GriesseMetlaRoesch09}).

Moreover, our results provided a little insight in the problem of
error estimates for joint Tikhonov-Lavrentiev regularization for
optimal control problems and obtained preliminary estimates. Further
research on the interpretation of the source condition and the Slater
condition in particular contexts seems necessary.  Finally we showed
that the main assumptions for our setting are fulfilled for a
distributed elliptic control problem and for some Fredholm equations
of the first kind.

\bibliographystyle{plain}
\bibliography{lorenz_roesch}
\end{document}